\documentclass[12pt]{amsart}

\usepackage[arrow,dvips,matrix,arrow,ps,color,line,curve,frame]{xy}

\usepackage{amssymb}

\advance\textheight1cm

\let\oldlabel=\label
\def\prellabel{\marginparsep=1em
    \def\label##1{\oldlabel{##1}\ifmmode\else\ifinner\else
         \marginpar{{\footnotesize\ \\ \tt
                    ##1}}\fi\fi}}

\advance\headheight1.15pt

\let\:\colon
\let\epsilon\varepsilon

\let\phi=\varphi
\let\theta=\vartheta

\let\Bbb=\mathbb
\let\frak=\mathfrak
\def\opn#1#2{\def#1{\operatorname{#2}}}
\opn\gp{gp} \opn\Max{Max} \opn\Ker{Ker} \opn\Coker{Coker}
\opn\Ext{Ext} \opn\conv{conv} \opn\chara{char} \opn\n{n} \opn\h{h}
\opn\GL{GuL} \opn\SL{SL} \opn\sn{sn} \opn\inte{int} \opn\End{End}
\opn\rank{rank} \opn\Aff{Aff} \opn\Spec{Spec} \opn\Proj{Proj}
\opn\QF{QF} \opn\I{Im} \opn\Hom{Hom} \opn\Aut{Aut} \opn\W{Witt}
\opn\W{W} \opn\inte{int} \opn\pyr{pyr} \opn\l{l} \opn\r{r}
\opn\const{const}

\def\ZZ{{\Bbb Z}}

\def\NN{{\Bbb N}}
\def\QQ{{\Bbb Q}}

\opn\End{End}
\opn\U{U}%
\opn\ch{ch}%

\def\Q{{\Box\kern1pt}}%

\def\c{\operatorname{\frak c}}
\def\kk{{\bf k}}

\newtheorem{lemma}{Lemma}

\newtheorem{theorem}[lemma]{Theorem}

\theoremstyle{definition}

\begin{document}

\title[Global coefficient ring in the nilpotence conjecture]{Global coefficient
ring in the\\ Nilpotence Conjecture}

\author{Joseph Gubeladze}

\thanks{Supported by NSF grant DMS-0600929}


\address{Department of Mathematics, San Francisco
State University, San Francisco, CA 94132, USA}

\email{soso@math.sfsu.edu}

\begin{abstract}
In this note we show that the nilpotence conjecture for toric
varieties is true over any regular coefficient ring containing
$\QQ$.
\end{abstract}

\maketitle


In \cite{G} we showed that for any additive submonoid $M$ of a
rational vector space with the trivial group of units and a field
$\kk$ with $\chara\kk=0$ the multiplicative monoid $\NN$ acts
nilpotently on the quotient $K_i(\kk[M])/K_i(\kk)$ of the $i$th
$K$-groups, $i\ge0$. In other words, for any sequence of natural
numbers $c_1,c_2,\ldots\ge2$ and any element $x\in K_i(\kk[M])$ we
have $(c_1\cdots c_j)_*(x)\in K_i(\kk)$ for all $j\gg0$
(potentially depending in $x$). Here $c_*$ refers to the group
endomorphism of $K_i(\kk[M])$ induced by the monoid endomorphism
$M\to M$, $m\mapsto m^c$, writing the monoid operation
multiplicatively.

The motivation of this result is that it includes the known
results on (stable) triviality of vector bundles on affine toric
varieties and higher $K$-homotopy invariance of affine spaces.
Here we show how the mentioned nilpotence extends to all regular
coefficient rings containing $\QQ$, thus providing the last
missing argument in the long project spread over many papers. See
the introduction of \cite{G} for more details.

Using Bloch-Stienstra's actions of the big Witt vectors on the
$NK_i$-groups \cite{St} (that has already played a crucial role in
\cite{G}, but in a different context), Lindel's technique of
\'etale neighborhoods \cite{L}, van der Kallen's \'etale
localization \cite{K}, and Popescu's desingularization \cite{Sw},
we show

\begin{theorem}\label{global}
Let $M$ be an additive submonoid of a $\QQ$-vector space with
trivial group of units. Then for any regular ring $R$ with
$\QQ\subset R$ the multiplicative monoid $\NN$ acts nilpotently on
$K_i(R[M])/K_i(R)$, $i\ge0$.
\end{theorem}

\medskip\noindent\emph{Conventions}. All our monoids and rings are
assumed to be commutative. $X$ is a variable. The monoid operation
is written mutliplicatively, denoting by $e$ the neutral element.
$\ZZ_+$ is the additive monoid of nonnegative integers. For a
sequence of natural numbers $\c=c_1,c_2,\ldots\ge2$ and an
additive submonoid $N$ of a rational space $V$ we put
$$
N^{\c}=\lim_{\to}\left(\xymatrix{N\ar[r]^{-^{c_1}}&N\ar[r]^{-^{c_2}}
&\cdots}\right)=
\bigcup_{j=1}^\infty N^{\frac1{c_1\cdots c_j}}\subset V.
$$

\

\begin{lemma}\label{graded}
Let $M$ be a finitely generated submonoid of a rational vector
space with the trivial group of units. Then $M$ embeds into a free
commutative monoid $\ZZ_+^r$.
\end{lemma}

For the stronger version of Lemma \ref{graded} with
$r=\dim_\QQ(\QQ\otimes M)$ see, for instance, \cite[Proposition
2.15(e)]{BG}. (In \cite{BG} the monoids as in Lemma \ref{graded}
are called the \emph{affine positive} monoids.)

\begin{lemma}\label{swanweibel}
Let $F$ be a functor from rings to abelian groups, $H$ be a monoid
with the trivial group of units, and $\Lambda=\bigoplus_H
\Lambda_h$ be an $H$-graded ring (i.~e.
$\Lambda_h\Lambda_{h'}\subset\Lambda_{hh'}$). Then we have the
implication
\begin{align*}
F(\Lambda)=F(\Lambda[H])\ \Longrightarrow\
F(\Lambda_e)=F(\Lambda).
\end{align*}
\end{lemma}

The special case of Lemma \ref{swanweibel} when $H=\ZZ_+$ is known
as the \emph{Swan-Weibel homotopy trick} and the proof of the
general cases makes no real difference, see \cite[Proposition
8.2]{G}.

\begin{lemma}\label{local}
Theorem \ref{global} is true for any coefficient ring of the form
$S^{-1}\kk[\ZZ_+^r]$ where $\kk$ is a field of characteristic $0$
and $S\subset\kk[\ZZ_+^r]$ is a multiplicative subset.
\end{lemma}

In the special case when $\kk$ is a number field Lemma \ref{local}
is proved in Step 2 in \cite[\S8]{G}, but word-by-word the same
argument goes through for a field $\kk$ provided the nilpotence
conjecture is true for the monoid rings with coefficients in
$\kk$.

\medskip\noindent\emph{Notice.} The reason we state the result in \cite{G} only
for number fields is that the preceding result in \cite{G} is the
validity of the nilpotence conjecture for such coefficient fields.
Actually, the proof of Theorem \ref{global} is an \'etale version
of the idea of interpreting the globalization problem for
coefficient rings in terms of the $K$-homotopy invariance, used
for Zariski topology in \cite[\S8]{G}.

\medskip Finally, in order to explain one formula we now summarize
very briefly the Bloch-Stienstra action of the ring of big Witt
vectors $\W(\Lambda)$ on
$$
NK_i(\Lambda)=\Coker(K_i(\Lambda)\to K_i(\Lambda[X])).
$$
For the details the reader is referred to \cite{St}.

The additive group of $\W(R)$ can be thought of as the
multiplicative group of formal power series $1+X\Lambda[[X]]$. It
has the decreasing filtration by the ideals
$I_p(R)=(1+X^{p+1}\Lambda[[X]])$, $p=1,2,\ldots,$ and every
element $\alpha(X)\in\W(\Lambda)$ admits a convergent series
expansion in the corresponding additive topology
$\alpha(X)=\Pi_{\NN}(1-\lambda_mX^m)$, $\lambda_m\in\Lambda$. To
define a continuous $\W(\Lambda)$-module structure on
$NK_i(\Lambda)$ it is enough to define the appropriate action of
the Witt vectors of type $1-\lambda X^m$, satisfying the condition
that every element of $NK_i(\Lambda)$ is annihilated by some ideal
$I_p(\W(\Lambda))$. Finally, such an action of $1-\lambda X^m$ on
$NK_i(\Lambda)$ is provided by the composite map in the upper row
of the following commutative diagram with exact vertical columns:
$$
\xymatrix{
0\ar[d]&0\ar[d]&0\ar[d]&0\ar[d]\\
NK_i(\Lambda[X])\ar[d]\ar[r]&NK_i(\Lambda[X])\ar[d]\ar[r]
&NK_i(\Lambda[X])\ar[d]\ar[r]&NK_i(\Lambda[X])\ar[d]\\
K_i(\Lambda[X])\ar[d]\ar[r]^{m^*}&K_i(\Lambda[X])\ar[d]\ar[r]^{\lambda_*}
&K_i(\Lambda[X])\ar[d]\ar[r]^{m_*}&K_i(\Lambda[X])\ar[d]\\
K_i(\Lambda)\ar[r]^{m\cdot-}\ar[d]&K_i(\Lambda)\ar[r]^{1}\ar[d]
&K_i(\Lambda)\ar[r]^{1}\ar[d]&K_i(\Lambda)\ar[d]\\
0&0&0&0}
$$
where:
\begin{itemize}
\item[(1)] $m_*$ corresponds to scalar extension through the
$\Lambda$-algebra endomorphism $\Lambda[X]\to\Lambda[X]$,
$X\mapsto X^m$,

\item[(2)] $m^*$ corresponds to \emph{scalar restriction} through
the same endomorphism $\Lambda[X]\to\Lambda[X]$,

\item[(3)] $\lambda_*$ corresponds to\ scalar extension through
the $\Lambda$-algebra endomorphism $\Lambda[X]\to\Lambda[X]$,
$X\mapsto\lambda X$.

\item[(4)] $m\cdot-$ is multiplication by $m$.
\end{itemize}

A straightforward check of the commutativity of the appropriate
diagrams, based on the description above, shows that for a ring
homomorphism $f:\Lambda_1\to\Lambda_2$ we have
\begin{equation}\label{witt}
f_*(\alpha z)=f_*(\alpha)f_*(z),\quad\alpha\in\W(\Lambda_1),\quad
z\in NK_i(\Lambda_1),
\end{equation}
where the same $f_*$ is used for the both induced homomorphisms
$$
\W(\Lambda_1)\to\W(\Lambda_2)\quad\text{and}\quad
NK_i(\Lambda_1)\to NK_i(\Lambda_2).
$$

\medskip\begin{proof}[Proof of Theorem \ref{global}] Since
$K$-groups commute with filtered colimits there is no loss of
generality in assuming that $M$ is finitely generated. Then by
Lemma \ref{graded} $R[M]$ admits a $\ZZ_+$-grading
$$
R[M]=\bigoplus_{\ZZ_+} R_j,\qquad R_e=R.
$$
In particular, by the Quillen local-global patching for higher
$K$-groups \cite{V}, we can without loss of generality assume that
$R$ is local.

\medskip\noindent\emph{Notice.} Actually, the local-global patching
proved in \cite{V} is for the special case of polynomial
extensions. However, the more general version for graded rings is
a straightforward consequence via the Swan-Weibel homotopy trick,
discussed above.

\medskip By Popescu's desingularization \cite{Sw} and the same
filtered colimit argument we can further assume that $R$ is a
regular localization of an affine $\kk$-algebra for a field $\kk$
with $\chara\kk=0$. In this situation Lindel has shown
\cite[Proposition 2]{L} that there is a subring $A\subset R$ of
the form $\kk[\ZZ_+^d]_\mu$, $\mu\in\max(\kk[\ZZ_+^d])$, $d=\dim
R$, such that
\begin{equation}\label{etale}
R\ \text{is \'etale over}\ A.
\end{equation}

\medskip\noindent\emph{Notice.} Lindel's result is valid in arbitrary
characteristic under the conditions that the residue field of $R$
is a simple separable extension of $\kk$, which is automatic in
our situation because $\chara\kk=0$.

\medskip Using again that $K$-groups commute with filtered colimits,
the validity of Theorem \ref{global} for $R$ is easily seen to be
equivalent to the equality
\begin{equation}\label{1}
K_i(R)=K_i(R[M^{\c}])
\end{equation}
for every sequence of natural numbers $\c=c_1,c_2,\ldots\ge2$.

Next we show that (\ref{1}) follows from the condition
\begin{equation}\label{2}
NK_i(R[M^{\c}])=0.
\end{equation}

In fact, by the filtered colimit argument we have
\begin{align*}
&K_i(R[M^{\c}])=K_i(R[M^{\c}])[X]\ \Longrightarrow\\
&K_i(R[M^{\c}])=K_i(R[M^{\c}])[\ZZ_+^{\c}]\
(=\lim_{\to}K_i(R[M^{\c}])[\ZZ_+]).
\end{align*}
On the other hand, by Lemma \ref{graded} the ring $R[M^{\c}]$ has
a $\ZZ_+^{\c}$-grading:
$$
R[M^{\c}]=\bigoplus_{\ZZ_+^{\c}}S_j,\qquad S_e=R.
$$
So by Lemma \ref{swanweibel} we have (\ref{1}).

\medskip To complete the proof it is enough to show (\ref{2}) assuming (\ref{etale}).

By the base change property the ring extension $A[M]\subset
R[M]=A[M]\otimes_A R$ is \'etale. Then by van der Kallen's result
\cite[Theorem 3.2]{K} we have the isomorphism of
$\W(R[M])$-modules
\begin{equation}\label{vdk}
NK_i(R[M])=\W(R[M])\otimes_{\W(A[M])}NK_i(A[M]).
\end{equation}

\medskip\noindent\emph{Notice.} Van der Kallen proves his formula
(the \emph{\'etale localization}) for a  modified tensor product
that takes care of the filtrations on $\W(A[M])$ and $\W(R[M])$ by
the ideals $I_p(-)$. But the presence of the characteristic $0$
subfield $\kk$ yields (via the \emph{ghost map}) the infinite
product presentations $\W(A[M])=\prod_{\NN}A[M]$ and
$\W(R[M])=\prod_{\NN}R[M]$ and, in particular, makes checkable the
appropriate compatibility of the two filtrations:
$I_p(R[M])=I_p(A[M])\W(R[M])$; see the discussion before
\cite[Theorem 3.2]{K}.

\medskip Pick an element $z\in NK_i(R[M])$. By (\ref{vdk}) it admits
a representation of the form
$$
z=\sum_q\alpha_q\bar y_q,\quad \alpha_q\in\W(R[M]),\quad y_q\in
NK_i(A[M]),
$$
where the bar refers to the image in $NK_i(R[M])$.

By Lemma \ref{local} we know that Theorem \ref{global} is true for
$A$. Therefore, $(c_1\cdots c_j)_*(y_q)=0$ for all $q$ provided
$j\gg0$. In particular, (\ref{witt}) implies
$$
(c_1\cdots c_j)_*(z)=\sum_q(c_1\cdots
c_j)_*(\alpha_q)\overline{(c_1\cdots c_j)_*(y_q)}=0,\qquad j\gg0.
$$
Since $z$ was an arbitrary element the filtered colimit argument
shows (\ref{2}).
\end{proof}

\end{document}